\documentclass[12pt]{amsart}
\usepackage{amsmath,amssymb,amsthm}
\newtheorem{theorem}{Theorem}
\newtheorem{lemma}[theorem]{Lemma}

\newtheorem{definition}{Definition}
\newcommand{\Ric}{\mathrm{Ric}}
\newcommand{\dvol}{d\mathrm{vol}}
\newcommand{\Vol}{\mathrm{Vol}}
\DeclareMathOperator{\conj}{conj}

\title[Conjugate radius and Ricci flow]{Conjugate radius and parabolicity along the Ricci flow}
\author{Chandan Kumar Mondal}
\address{Department of Mathematics, School of Sciences, Netaji Subhas Open University, West Bengal, India}
\email{chan.alge@gmail.com}
\begin{document}

\begin{abstract}
We study the behaviour of the conjugate radius under the Ricci flow on complete noncompact manifolds. We show that if $(M^n,g(t))_{t\in[0,T]}$ is a
complete Ricci flow with nonnegative scalar curvature at $t=0$, the initial
metric is parabolic and not scalar-flat, and the Einstein tensor
$\mathrm{Ric}-\tfrac12 Rg$ is bounded above, then parabolicity is preserved
along the flow and the conjugate radius of $g(t)$ is finite for every
$t\in[0,T]$. As an application, we prove that a
complete noncompact gradient shrinking Ricci soliton whose volume growth is subquadratic, in the sense that $\liminf_{r\to\infty}\Vol(B(p,r))/r^2=0$, has finite conjugate radius.
\end{abstract}
\noindent\footnotetext{$\mathbf{2020}$\hspace{5pt}Mathematics\; Subject\; Classification: Primary 53C20, 53C21, 53E20. }

\maketitle

\section{Introduction}

In 1982, Hamilton \cite{HA82} introduced the Ricci flow, which later became a central tool in the program toward the Poincar\'e conjecture. The Ricci flow is defined by an evolution equation for metrics on a Riemannian manifold $(M,g_0)$:
\begin{equation}\label{r1}
   \frac{\partial}{\partial t}\,g(t)=-2\,\mathrm{Ric}_{g(t)},\qquad g(0)=g_0 .
\end{equation}
Hamilton proved short-time existence of the Ricci flow on closed manifolds. For complete metrics with bounded curvature, short-time existence was established by Shi \cite{Shi89}, and uniqueness of solutions was proved by Chen and Zhu \cite{CZ06} under the bounded curvature condition.

\par The conjugate radius of a manifold is the largest $r$ such that the exponential map is nonsingular on every ball of radius $r$ in the tangent spaces. Estimating the conjugate radius is important, since within this radius geodesics have no conjugate points. A complete Riemannian manifold $(M,g)$ has no conjugate points precisely when its conjugate radius $\conj(M,g)$ is infinite. Sectional curvature plays an important role in estimating the conjugate radius; for example, the Cartan--Hadamard theorem states that a complete manifold with nonpositive sectional curvature has infinite conjugate radius. Besides sectional curvature bounds, scalar curvature bounds also affect the estimation of the conjugate radius. One important result is due to Green \cite{Gr58}, who gave a bound on the conjugate radius for closed manifolds with positive average scalar curvature. Zhu \cite{Zhu22} extended this result to the case of open manifolds with an additional nonnegative Ricci curvature assumption.
Very recently, Kwong \cite{Kwo25} and Ge and Zhang \cite{GZ26} have estimated the conjugate radius under various scalar curvature assumptions.

In the case of the Ricci flow, the conjugate radius is not always preserved, even on closed manifolds. For example, on the double torus $\Sigma_2$ ($\chi=-2$) we can choose a metric $g_0$ with finite conjugate radius. Let $g(t)$ be the normalized Ricci flow with $g(0)=g_0$. Then $g(t)$ converges to a smooth metric $g_\infty$ of constant curvature $k$, conformal to $g_0$; see \cite{HAM88}. Since $\chi(\Sigma_2)=-2$, the Gauss--Bonnet theorem implies that $k<0$. But by the Cartan--Hadamard theorem, $(\Sigma_2, g_\infty)$ has infinite conjugate radius.

Therefore, a suitable curvature condition is necessary for the Ricci flow to preserve finiteness of the conjugate radius. Motivated by these observations, we study how the conjugate radius behaves under the Ricci flow on complete noncompact manifolds. In this article, we assume the existence of the Ricci flow on noncompact manifolds. The following is the main theorem of this article:

\begin{theorem}\label{thm:C}
Let $(M^n,g(t))_{t\in[0,T]}$, with $n\ge2$ and $T>0$, be a solution of the Ricci flow on a connected noncompact manifold $M$, with every $g(t)$ complete, such that $R(\cdot,0)\ge0$, and $R(\cdot,0)\not\equiv 0$. Assume that $(M,g(0))$ is parabolic and that for some $K_0\in\mathbb R$
\[
G_{g(t)}=\Ric_{g(t)}-\tfrac12R_{g(t)}\,g(t)\le K_0\,g(t)\qquad\text{on }M\times[0,T].
\]
Then for every $t\in[0,T]$:
\begin{enumerate}
\item[(i)] $(M,g(t))$ is parabolic;
\item[(ii)] $\conj(M,g(t))<\infty$.
\end{enumerate}
\end{theorem}

We do not assume that the curvature is bounded, nor that the flow is unique. The only curvature restrictions we need are nonnegativity of the scalar curvature of the initial metric and an upper bound on the Einstein tensor along the flow. This result also shows that parabolicity is preserved along the Ricci flow under suitable conditions. In Section 2, we give the proof of the main theorem. In Section 3, we show that a complete noncompact gradient shrinking Ricci soliton with subquadratic volume growth has finite conjugate radius.
\subsection*{AI disclosure}
The author used ChatGPT as an assistant during the preparation of this manuscript to edit \LaTeX\ text, check calculations, and refine the manuscript language. All
proofs were written entirely by the author, who retains exclusive responsibility for the content of this work.

\section{Proof of the main result}\label{sec:parab}

Let $\mathrm{Lip}_c(M)$ denote the compactly supported locally Lipschitz functions on $M$. For a metric $g$ and a compact set $K\subset M$, the capacity is defined by, see \cite{GS02},
\[
\operatorname{cap}_g(K)=\inf\Big\{\int_M|\nabla f|^2_g\,\dvol_g:\ f\in\mathrm{Lip}_c(M),\ f\ge1\text{ on }K\Big\}.
\]
\begin{definition}\cite{GS02, Gri99}
    A complete noncompact manifold $(M,g)$ is \emph{parabolic} if $\operatorname{cap}_g(K)=0$ for every compact $K\subset M$. Equivalently, every smooth positive superharmonic function on $M$
(i.e.\ $u>0$, $\Delta_g u\le 0$) is constant. Otherwise, $(M,g)$ is called
\emph{nonparabolic}.
\end{definition}

To prove the main result, we need the inequality of Lemma 3.1 of \cite{GZ26}. 

\begin{lemma}\cite{GZ26}\label{lem:ineq}
Let $(M^n,g)$ be complete with $n\ge2$, and let $0<a\le\conj(M,g)$ be finite. Then for every compactly supported Lipschitz function $f$ on $M$,
\begin{equation}\label{eq:1}
\int_M R\,f^2\,\dvol\;\le\; n(n-1)\frac{\pi^2}{a^2}\int_M f^2\,\dvol+(n-1)\int_M|\nabla f|^2\,\dvol .
\end{equation}
\end{lemma}
Next we show that parabolicity is preserved by the Ricci flow, as long as the Einstein tensor $G=\Ric-\tfrac12Rg$ stays bounded from above. 

\begin{proof}[Proof of Theorem \ref{thm:C}]
Since for any fixed $K>0$, we have $R(\cdot,0)\ge 0\ge -K$, the Corollary 2.3(i) of \cite{CH09} gives
$R(\cdot,t)\ge -\dfrac{n}{2t+n/K}$ on $M\times[0,T]$, because $g(t)$ is a complete
smooth solution. Letting $K\to 0$ yields $R\ge 0$ on $M\times[0,T]$.
Moreover, the evolution equation of scalar curvature under Ricci flow is $$\partial_t R=\Delta R+2|\mathrm{Ric}|^2\ge \Delta R.$$ Hence, $R\ge 0$ is a
supersolution of the heat equation. If $R(p,t_0)=0$ for some $p\in M$ and
$t_0\in(0,T]$, the strong maximum principle \cite{Nir53}, applied in local coordinate charts
and combined with the connectedness of $M$, gives $R\equiv 0$ on $M\times[0,t_0]$, contradicting
$R(\cdot,0)\not\equiv 0$. Hence $R(\cdot,t)>0$ on $M$ for all $t\in(0,T]$.

To prove the first part, consider the energy function
\begin{equation}\label{g1}
E_t(f)=\int_M|\nabla f|^2_{g(t)}\,\dvol_{g(t)},
\end{equation}
for $f\in\mathrm{Lip}_c(M)$. Under the Ricci flow, we have the following evolution equations; see \cite{HA82}:
\[
\partial_tg^{ij}=2R^{ij},\qquad \partial_t\,\dvol_{g(t)}=-R\,\dvol_{g(t)}.
\]
Since a Lipschitz function is differentiable almost everywhere, at every point where $df$ exists, we calculate
\begin{align}\label{eq:integrand}
\nonumber \partial_t\Big(|df|^2_{g(t)}\,\dvol_{g(t)}\Big)
&=\big(\partial_tg^{ij}\big)\,\partial_if\,\partial_jf\,\dvol_{g(t)}+g^{ij}\,\partial_if\,\partial_jf\,\big(\partial_t\dvol_{g(t)}\big)\\
\nonumber &=2R^{ij}\,\partial_if\,\partial_jf\,\dvol_{g(t)}-R\,g^{ij}\,\partial_if\,\partial_jf\,\dvol_{g(t)}\\
\nonumber &=\big(2R^{ij}-R\,g^{ij}\big)\,\partial_if\,\partial_jf\,\dvol_{g(t)}\\
&=2\,G_{g(t)}(\nabla f,\nabla f)\,\dvol_{g(t)},
\end{align}
where we have used that $G^{ij}=R^{ij}-\tfrac12Rg^{ij}$.
\par Let $S=\operatorname{supp}f$. Then $S$ is a compact set. By Rademacher's theorem \cite[\S5.8.3]{EL10}, $df$ exists almost everywhere and is bounded on $S$. Write $\dvol_{g(t)}=\rho(\cdot,t)\,\dvol_{g(0)}$, where $\rho$ is smooth and positive and satisfies $\partial_t\rho=-R\rho$. On the compact set $S\times[0,T]$ the metrics $g(t)$, their inverses, the function $\rho$ and the curvature are all bounded with respect to $g(0)$. By \eqref{eq:integrand} there is a constant $C$ with
\[
\Big|\partial_t\big(|df|^2_{g(t)}\,\rho(\cdot,t)\big)\Big|=2\Big|\,G_{g(t)}(\nabla f,\nabla f)\,\rho(\cdot,t)\Big| \le C\,|df|^2_{g(0)}\le C',
\]
and the bound $C'$ is independent of $t$ and integrable over $S$. By dominated convergence, passing the limit inside the integral, we get
\[
\frac{d}{dt}E_t(f)=2\int_MG_{g(t)}(\nabla f,\nabla f)\,\dvol_{g(t)} .
\]
Since $G_{g(t)}(\nabla f,\nabla f)\le K_0|\nabla f|^2_{g(t)}$ almost everywhere, this gives $\frac{d}{dt}E_t(f)\le2K_0\,E_t(f)$, hence
\[
\frac{d}{dt}\Big(e^{-2K_0t}E_t(f)\Big)=e^{-2K_0t}\Big(\frac{d}{dt}E_t(f)-2K_0E_t(f)\Big)\le0 .
\]
So $e^{-2K_0t}E_t(f)$ is nonincreasing. This yields
\begin{equation}\label{eq:energy}
E_t(f)\le e^{2K_0t}E_0(f)\qquad \text{for }t\in[0,T].
\end{equation}
Now let $K\subset M$ be compact and $\varepsilon>0$. Choose $f\in\mathrm{Lip}_c(M)$ with $f\ge1$ on $K$ and $E_0(f)\le\operatorname{cap}_{g(0)}(K)+\varepsilon$. The same $f$ is admissible for $\operatorname{cap}_{g(t)}(K)$, since the admissible class involves neither the metric nor the time. By \eqref{eq:energy},
\[
\operatorname{cap}_{g(t)}(K)\le E_t(f)\le e^{2K_0t}\big(\operatorname{cap}_{g(0)}(K)+\varepsilon\big).
\]
Letting $\varepsilon\to0$ and using that $(M,g(0))$ is parabolic, i.e.\ $\operatorname{cap}_{g(0)}(K)=0$, we get $\operatorname{cap}_{g(t)}(K)=0$ for every compact $K$. Since each $g(t)$ is complete, $(M,g(t))$ is parabolic, which proves (i).

To prove the second part, fix $t\in[0,T]$ and suppose $\conj(M,g(t))=\infty$. Then Lemma~\ref{lem:ineq} applies to $g(t)$ with every finite $a>0$. For $f\in\mathrm{Lip}_c(M)$ the integral $\int_Mf^2\,\dvol$ is finite, so letting $a\to\infty$ yields
\begin{equation}\label{eq:Rf}
\int_M R\,f^2\,\dvol_{g(t)}\le(n-1)\int_M|\nabla f|^2_{g(t)}\,\dvol_{g(t)},
\end{equation}
for all $f\in\mathrm{Lip}_c(M)$.

Pick $x\in M$ such that $R(x,0)>0$ and let $K=\overline{B}_{g(t)}(x,1)$, which is compact because $g(t)$ is complete. Given $\varepsilon>0$, part (i) provides $f\in\mathrm{Lip}_c(M)$ with $f\ge1$ on $K$ and $\int_M|\nabla f|^2_{g(t)}\,\dvol_{g(t)}<\varepsilon$. Because $R\ge0$ and $f^2\ge1$ on $K$, inequality \eqref{eq:Rf} gives
\[
\int_K R\,\dvol_{g(t)}\le\int_M R\,f^2\,\dvol_{g(t)}\le(n-1)\,\varepsilon .
\]
Letting $\varepsilon\to0$ we get $\int_K R\,\dvol_{g(t)}=0$. On the other hand, $K$ contains an open ball with $R(\cdot,t)>0$, so $\int_KR\,\dvol_{g(t)}>0$. This contradiction shows that $\conj(M,g(t))<\infty$.
\end{proof}

\section{Gradient shrinking Ricci solitons}\label{sec:soliton}
A self-similar solution of the Ricci flow is called a Ricci soliton. More precisely, a Ricci flow $g(t)$ on $M$ is a Ricci soliton if
$$g(t)=\sigma(t)\,\varphi_t^*\,g_0$$
for a positive function $\sigma(t)$ with $\sigma(0)=1$ and diffeomorphisms $\varphi_t$ with $\varphi_0=\mathrm{id}$. Differentiating the above equation, we obtain an equation of the form
$$\mathrm{Ric}_{g_0}+\tfrac12\mathcal L_Xg_0=\lambda g_0$$
for a vector field $X$ and a constant $\lambda$. When $X$ is the gradient vector field of a smooth function $f\in C^\infty(M)$, the above equation takes the form
\begin{equation}\label{r2}
 \mathrm{Ric}+\nabla^2f=\lambda g.
\end{equation}
When a manifold $(M,g)$ satisfies equation (\ref{r2}), it is called a gradient Ricci soliton. The Ricci
soliton is called shrinking, steady, or expanding according as $\lambda>0$, $\lambda = 0$, or $\lambda<0$, respectively.
\par Hamilton showed (see \cite{ELM08}) that the following identity holds for a gradient Ricci soliton:
\begin{equation}\label{eq:solid}
\Delta R-\langle\nabla f,\nabla R\rangle=2\lambda R-2|\Ric|^2 .
\end{equation}
The following lemma is needed to prove the main result of this section.

\begin{lemma}\label{th1}
Let $(M^n,g)$ be a complete Riemannian manifold with $n\ge2$, scalar curvature $R\ge0$ and $\conj(M,g)=\infty$. Suppose that for some $p\in M$
\begin{equation}\label{eq:growth}
\liminf_{r\to\infty}\frac{\Vol(B(p,r))}{r^2}=0 .
\end{equation}
Then $R\equiv0$.
\end{lemma}

\begin{proof}
For a fixed $r>0$, consider the Lipschitz cut-off function $\varphi_r$ defined in \cite{Yau76} such that
\[
\begin{cases}
0\le\varphi_r\le1 & \text{on } M,\\[2pt]
\varphi_r=1 & \text{in } B(p,r),\\[2pt]
|\nabla\varphi_r|\le\dfrac{1}{r} & \text{almost everywhere},\\[6pt]
\operatorname{supp}(\varphi_r)\subset B(p,2r).
\end{cases}
\]

The gradient of a Lipschitz function vanishes almost everywhere on its level sets, so $\nabla\varphi_r$ vanishes almost everywhere outside $\{r<d(p,\cdot)<2r\}\subset B(p,2r)$. Therefore, we get
\begin{eqnarray}\label{eq:cutoff}
\nonumber \int_M\varphi_r^2\,\dvol &\le &\Vol(B(p,2r)),\qquad \text{and}\\
\int_M|\nabla\varphi_r|^2\,\dvol & \le &\frac{1}{r^2}\Vol(B(p,2r)).
\end{eqnarray}
Since $\conj(M,g)=\infty$, Lemma~\ref{lem:ineq} holds for every finite $a>0$. Applying it with $f=\varphi_r$ and $a=r$, and using \eqref{eq:cutoff}, we get
\[
\int_M R\,\varphi_r^2\,\dvol\le\frac{c_n}{r^2}\,\Vol(B(p,2r)),\qquad c_n=n(n-1)\pi^2+(n-1).
\]
Because $R\ge0$ and $\varphi_r\equiv1$ on $B(p,r)$, for every $r>0$ we have
\begin{equation}\label{eq:step2}
\int_{B(p,r)}R\,\dvol\le\int_MR\,\varphi_r^2\,\dvol\le\frac{c_n}{r^2}\,\Vol(B(p,2r)).
\end{equation}
By \eqref{eq:growth} there is a sequence $s_k\to\infty$ with $\Vol(B(p,s_k))/s_k^2\to0$. Taking $r_k=s_k/2$ in \eqref{eq:step2},
\[
\int_{B(p,r_k)}R\,\dvol\le 4c_n\,\frac{\Vol(B(p,s_k))}{s_k^2}\longrightarrow0 .
\]
Since $R\ge0$, the integral $\int_{B(p,r)}R\,\dvol$ is nondecreasing in $r$ and converges to $\int_MR\,\dvol$ by monotone convergence. Hence $\int_MR\,\dvol=0$, and as $R\ge0$ is continuous, $R\equiv0$.
\end{proof}

\begin{theorem}\label{prop:soliton}
Let $(M^n,g,f)$, $n\ge2$, be a complete noncompact gradient shrinking Ricci soliton. If $\liminf_{r\to\infty}\Vol(B(p,r))/r^2=0$ for some $p\in M$, then $\conj(M,g)<\infty$.
\end{theorem}

\begin{proof}
By Chen \cite[Cor.~2.5]{CH09}, a complete gradient shrinking soliton has $R\ge0$. Suppose that $\conj(M,g)=\infty$. Then Lemma~\ref{th1} gives $R\equiv0$, and \eqref{eq:solid} yields $|\Ric|^2\equiv0$. Therefore \eqref{r2} implies $\nabla^2f=\lambda g$ with $\lambda>0$. Since $f$ is nonconstant, Tashiro's theorem \cite{Ta65} implies that $M$ is isometric to $\mathbb R^n$. But in $\mathbb R^n$ we have $\Vol(B(p,r))/r^2=\omega_nr^{n-2}\ge\omega_n>0$ for $r\ge1$, where $\omega_n$ is the volume of the unit ball, which contradicts the volume assumption. Hence $\conj(M,g)<\infty$.
\end{proof}

\end{document}